\documentclass[12pt]{mycls}

\usepackage{graphicx}
\usepackage{amsmath}
\usepackage{amssymb}
\usepackage{latexsym}
\usepackage{subfigure}

\newcommand{\be}{\begin{equation}} \newcommand{\ee}{\end{equation}}
\newcommand{\barr}{\begin{array}} \newcommand{\earr}{\end{array}}
\newcommand{\bitem}{\begin{itemize}}
\newcommand{\eitem}{\end{itemize}}

\newcommand{\dstyle}{\displaystyle}

\newtheorem{thm}{Theorem}[section]
\newtheorem{lemma}[thm]{Lemma}
\newtheorem{alg}{Algorithm} 

\begin{document}
\begin{frontmatter}

\title{A Non-Krylov subspace Method for Solving Large and Sparse Linear System of Equations}
\author[au1]{Wujian Peng\corref{cor1}}    \ead{wpeng@zqu.edu.cn}
\author[au2]{Qun Lin} \ead{linq@lsec.ac.cc.cn}
\address[au1]{Department of Math. and stats. Sciences,
          Zhaoqing Univ.,  Zhaoqing,  China,526061 }
\address[au2]{Academy of Math. and System Sciences,
           Chinese Academy of Sciences, China,100081}
\cortext[cor1]{Corresponding Author}

\begin{abstract}
   Most current prevalent iterative methods can be classified into the so-called extended Krylov subspace methods, a class of iterative methods which do not fall into this category are also proposed in this paper. Comparing with traditional Krylov subspace methods which always depend on the matrix-vector multiplication with a fixed matrix, the newly introduced methods(the so-called  (progressively) accumulated projection methods, or AP (PAP) for short) use a projection matrix which varies in every iteration to form a subspace from which an approximate solution is sought. More importantly an accelerative approach(called APAP) is introduced to improve the convergence of PAP method. Numerical experiments demonstrate some surprisingly improved convergence behavior. Comparison between benchmark extended Krylov subspace methods(Block Jacobi and GMRES) are made and one can also see remarkable advantage of APAP in some examples. APAP is also used to solve systems with extremely ill-conditioned coefficient matrix (the Hilbert matrix) and numerical experiments shows that it can bring very satisfactory results even when the size of system is up to a few thousands.
\end{abstract}

\end{frontmatter}

\section{Introduction}

  Linear systems of the form
  \be\label{eq:1}
  Ax =b
  \ee where $A\in R^{n \times n}$ being nonsingular arise from tremendous mathematical applications and are the
  fundamental objects of almost every \mbox{computational} process. From the very ancient
  Gaussian elimination to the state-of-the-art methods like CG, MINRES, GMRES,
  as well as Multigrid method\cite{Axelsson1, Axelsson2, templates, Golub, Hackbusch_MG,Hackbusch_It}, numerous solvers of
  linear systems have been introduced and studied in
  extreme detail. Basically all solvers fall into two categories: direct
  methods  and   iterative methods.

  Except for those specially designed methods for systems with some special
  properties, like symmetry, sparsity or triangularity, elimination
  methods based on LU factorization seem to be most widely accepted for  general
  linear systems with satisfactory stability due to its flexibility of pivoting  strategies\cite{SuperLU_smp99, Duff1984, Alan}. Comparing with  direct methods, iterative methods are a much larger family and have been accepting dominant attention. Since they make it possible for people to get a very `close' solution to a system in much less arithmetic operation and storage requirement than direct methods and thus often lead to huge savings of time and costs.

  Although some state-of-the-art direct methods can be applied to solve systems with pretty large amount of unknowns\cite{templates, Duff1997} in some situations, for even larger scale sparse systems(say, with unknowns up to a few millions) one can resort to the LGO-based solver\cite{pengXueBao, pengDDM2009} recently introduced by authors, iterative methods are the  only option available for many practical problems. For example, detailed three-dimensional multiphysics
  simulations lead to linear systems comprising hundreds of millions or even billions of equations in as many
  unknowns, systems with several millions of unknowns are now routinely encountered in many applications, making  the use of iterative methods virtually mandatory.

  The history of iterative methods can largely be divided into two major periods. The first period begins with 1850's  while Jacobi and Gauss etc. established the first iterative methods named after these outstanding researchers and the period ends in 1970's. The majority of these iterative method are classified as stationary methods, which usually take the form:
      \be \label{eq:stationaryIt}
           x_{k+1} = G x_k + v,  \quad (k = 0,1,2, \cdots).
      \ee
 where $v$ is a fixed vector and $x_0$ as the first guess.
  Excellent books covering the detailed analysis of error and convergence of these methods include works by Axelsson\cite{Axelsson2}, Datta\cite{Datta1995book}, Varga\cite{Varga}
  and David Young\cite{David}, etc.
  The second period begins in the mid-1970s and is dominated by Krylov subspace methods and preconditioning techniques. Generally Krylov subspace methods use the following form
   \be
       x_{k} = x_0 + y_k, \quad (k =1,2,\cdots)
   \ee
where  $x^0$ is an initial guess and $y^k$ belongs to a so-called Krylov subspace
$$ {\cal K}_m(G,v_0) \equiv span\{ v_0, Gv_0, G^2 v_0, \cdots, G^{m-1}v_0  \}. $$
By assuming different strategies for seeking $y_k$ from $ {\cal K}_m(G,v_0)$, one gets a variety of
iterative methods such as CG, BiCG, GMRES, FOM, MINRES, SYMMLQ, QMR\cite{Freund1991,paige,Saad,saad1986,Vorst1992}, etc.

As a matter of fact, if we would refer extended Krylov subspace methods as those at each step of iteration the correction vector or approximate solution always comes from Krylov subspaces with a few fixed ``generator" matrices (by a ``generator" matrix to Krylov subspace ${\cal K}_m(A,v)$ we mean matrix $A$ here), then the traditional stationary iterative methods such as Jacobi, Gauss-Seidal, SOR as well as the more general Richardson iterative methods can also be classified as extended Krylov subspace methods. Since for example one can easily see from (\ref{eq:stationaryIt}) that
$$
\begin{array}{rl}
x_{k+1}  &= v + Gx_k = v+ G( v+Gx_{k-1} ) \\
    & = v + Gv + G^2x_{k-1}=v + Gv + G^2(v + Gx_{k-2}) \\
    &\cdots\\
    &= v + Gv + G^2v + G^3 v + \cdots + G^{k}v + G^{k+1 }x_0 \equiv y_{k+1} + z_{k+1}
\end{array}
$$
where $y_{k+1} = v + Gv + G^2v + G^3 v + \cdots + G^{k}v \in {\cal K}_{k+1}(G,v)$ and $z_{k+1}=G^{k+1}x_0 \in {\cal K}_{k+2}(G,x_0)$ and $x_0$ is the initial guess to the system. In a word, any iterative scheme that takes the following form
\be\label{iter:general_form}
x_{k+1} = \sum_1^m\mathcal{P}_i(G)v_i
\ee
can be classified into the extended Krylov subspace methods, where $\mathcal{P}_i(G)$ ($i=1,2,\cdots, m$) denotes a matrix polynomial function, $G$ is the so-called iterative matrix and $v_i$ ( $i=1,2,\cdots, m$) is usually some fixed starting vector, and $m$( usually two) is a very small integer.

We need to mention that the well-known row projection methods such as Karczmarz's method(known as ART method in computed tomography) and Cimmino's methods can also be regarded as stationary iterative methods\cite{Bramley,Galantai}, thus they also belong to the category of extended Krylov subspace methods.

The extended Krylov subspace methods may  be very effective when the coefficient matrix is close to the normal matrix,  or the exact solution lie on Krylov subspace formed by the eigenvectors corresponding the leading eigenvalues in magnitude.  However since the base vectors of Krylov subspaces always take the form $G^k v$, it can be very inefficient to find a good ``approximation" to the error vector in such a subspace when the condition number of the coefficient matrix is large, especially when vector is almost perpendicular to the Krylov subspace ${\cal K}_m(A,v)$.  Thus even when extended Krylov subspace methods are applied on relatively small-sized systems, in many cases one still has to use some kind of preconditioning techniques to obtain an improved convergence.

It is always desirable for us to use some types of preconditioning techniques when we apply iterative methods to solve linear system of \mbox{equations}, especially for large scale computing. Though  numerous preconditioning techniques are exploited in recent decades and some of them turn out to be extremely efficient in some special situations, there does not exist a simple preconditioning technique which can be applied in  general cases. Another important factor is, all preconditioning techniques can be traced back to certain algebraic iterative schemes\cite{Benzi, Vorst}.

It is therefore our motivation here to develop a set of purely
algebraic algorithms that can in someway overcome the difficulties
arising in the extended Krylov subspace methods. In the meantime we also develop
some accelerating techniques to improve the convergence of our new
iterative methods. Our intention here is, instead of using Krylov
subspace methods with a fixed generator matrix and fixed starting
vector we use a sequence of subspaces formed by some base vectors
that are eventually approximating the exact solutions. Since the
base vectors are obtained by some successive projections and has the
property that it carries the largest magnitude in some subspaces, we
name them as Accumulated Projection Methods(AP).

\section{Basic Principles for Iterative Methods}
In this section we review the basic rules that govern the designing of iterative methods for solving linear system equations,
which in turn helps to derive our methods introduced in later sections.

Currently any iterative solver for system (\ref{eq:1}) always begins with an initial guess $x_0$
(without assumptions imposed on $x_0$), which leaves an easily available residual vector $r_0$ defined as $r_0 = b - A x_0$.
If we denote the error vector as $e_0 = x - x_0$, we then have $ Ae_0 = r_0$.  An effective iterative scheme then
seeks a sequence of vector $\{x_k\}_0^\infty$ so that the corresponding sequence $\{e_k\}_0^\infty $
(with $e_k = x - x_k, \, k = 0,1,2\cdots, \infty $) will converge to zero vector in $R^n$, or equivalently
the sequence of error norms $\{||e_k||\}_0^\infty$ converge to zero. If we assume that the coefficient
matrix $A$ in system (\ref{eq:stationaryIt}) is nonsingular, one can see that the sequence of
residual norms $\{||r_k||\}_0^\infty $ (with $r_k = b - Ax_k, k = 0,1,2\cdots, \infty $) also converges to
zero since we always have  $ Ae_k = r_k$ for $ k=0,1,2,\cdots$, which leads to $ ||r_k|| \le ||A|| ||e_k||$.
For example, the traditional stationary iterative methods such as Jacobi, Gauss-Seidal and SOR methods satisfy (\ref{eq:stationaryIt})
with the iterative matrix $G$ taken different form in each situation, and to make these iterative scheme convergent,
a sufficient and necessary condition is
$$ \rho(G)< 1.$$
Note that from (\ref{eq:stationaryIt}) we have $e_{k+1} = G e_k$, this implies that the sequence of error norms
 is strictly decreasing and has zero as its limit. In Krylov subspace methods, people usually expect either the sequence of error norms(in CG, this is the $||\cdot||_a$ of the error defined by $||e||_a = ||e^T A e||$\cite{templates})  or residual norms $ ||r||$( In GMRES, this is the regular $ ||\cdot ||_2$ norm) are decreasing sequences and converge to zero.  It should be kept in mind in general a small value of residual norm can not be used as an indication of convergence for an iterative process, while direct estimation of error norms is practically not possible, thus one often uses the relative residual norm as its convergence indicator.

Traditional iterative schemes of the form (\ref{eq:stationaryIt}) usually depend on the splitting of coefficient matrix $A$, while effective ways of splitting of $A$ which lead to convergent iterative schemes usually require $A$ satisfying certain special property(diagonally dominant, SPD, etc.) and thus not so easy to design. Many well-established  iterative schemes(including CG, MINRES, SYMMLQ) need special properties of $A$(SPD, or symmetry, etc.); only a few well-known iterative methods(GMRES, BiCG, LSQR etc.) can be applied to general nonsingular coefficient matrices and unfortunately none of these methods have well-established convergence analysis. Since linear systems of \mbox{equations} come from various scientific computation and engineering practicing, the required properties for many of these iterative schemes can not be satisfied in general, it is thus more attractive to design iterative methods for general linear system of equations.

In the later sections, we will apply the basic principles to design
a convergent iterative scheme for solving system (\ref{eq:1}),
specifically we will use projection techniques to get a sequence of
approximations $\{x_k\}$ to exact solution $x$ so that the error
vectors  $e_k$ ($e_k=x-x_k$) have strictly decreasing Euclidean
norms. We will use a strategy which differs from any current Krylov
subspace methods. First of all in our method the initial guess
vector $x_0$ to the solution of (\ref{eq:1}) can not be chosen
arbitrarily, instead we suggest a few ways to construct a ``good"
initial guess, in later searching of corrections to previous approximations we
don't use any Krylov subspaces and there is no so-called iterative matrix (like $G$ in (\ref{eq:stationaryIt}))
in the whole process, thus they do not fall into the category of extended
Krylov subspace methods.

\section{An Accumulated Projection Idea}

In essence  every iterative scheme always tries to seek an approximate solution in as less as possible steps.
Equivalently we wish to construct a subspace with much smaller dimension than $n$( the number of unknowns in the
system) and then seek a good approximate solution in this subspace. Currently all prevalent iterative schemes use
one or two fixed generator matrices to create one or two Krylov subspaces frow where an approximated solution(correction) may be obtained in these subspaces. However in practical computation  Krylov subspace ${\cal K}_m(A,v)$ always stays close to the leading eigenspace $L_s(A)$ defined by
$$ L_s(A)=span(v_1,v_2,\cdots, v_s)$$ where $Av_i = \lambda_i v_i$ and $\lambda_i(i=1,2,\cdots s)$ are the largest eigenvalues of $A$ in terms of magnitude in decreasing order, since in finite precision computing various computing errors(rounding-off errors, errors caused by cancellation of significant digits, etc) can not be avoided, especially in large scale computation. We will present a different approach to construct a subspace where no adoption of any vectors in the form $G^k v$ for its basis vectors is used and thus we can expect to avoid the drawbacks related to this type of subspaces.

Let's start from a simple projection idea. If we check each of the
row in system (\ref{eq:1}) we have $a_i x = b_i$, where $a_i$ is the
$i$-th row vector of the coefficient matrix $A$ and $b_i$ is the
$i$-th component of the right side vector $b$. A natural idea is to use
the projection vector $p_i$ of $x$ ($p_i = \alpha a_i^T$) on the
direction $a_i$ as its approximation, where $\alpha = \frac{b_i}{
a_ia_i^T}$. The corresponding error vector $e_i = x - p_i  $
satisfies \be  \label{errorNorm} ||e_i||^2 = ||x||^2 - ||p_i ||^2.
\ee
A simple successive application of this process gives the so-called
Row Projection Methods first proposed by Karcmarz, and was later found that they
are nothing but a stationary iterative method:

\be \label{def:RP} x_{k+1} = Q_u x_k + b_u \ee
where the iterative matrix $Q_u$ is formed as $$ Q_u = (I -  P_m)(I -  P_{m- 1}) \cdots(I -  P_1).$$
and $P_i$ ($ 1\le i \le m$) are projection vectors to some subspaces of $R^n$.
Another type of  Row projection approach for solving (\ref{eq:1}) is proposed by Cimmino in 1939\cite{cimmino}. Cimmino's approach was later found to be
equivalent as block Jacobian iteration  with the iterative matrix having the form
$$ G = \sum_i \omega_i P_i$$ where $P_i$ represents the projection matrix over some subspaces formed by some row vectors of matrix $A$ and $\omega_i$ are some carefully chosen parameters so that $\rho(G)<1$.
These row projection methods have been examined by several authors and some accelerative schemes are proposed to improve the convergence behavior\cite{Bramley,Galantai}.

In the following subsections we are to present a new type of projection technique---accumulated projection. Unlike the row projection techniques which end up with the form of some stationary
iterative schemes\cite{Galantai} and thus fall into the extended Krylov subspace methods,  our AP technique does not depend on any Krylov subspace.

\subsection{An accumulated projection}
The best approximation vector to $x$ in terms of error length(i.e., its Euclidean norms) in any subspace $W$ of $R^n$ is its projection $p (\in W)$. In exact arithmetic, the bigger the dimension of $W$ is, the bigger the length of $p$, i.e.,
the closer the two vectors $x$ and $p$ in terms of their angle.  Unfortunately in practical computation if $W$ is usually constructed by using Krylov subspace technique with a fixed generator matrix $A$, i.e., $W= {\cal K}_m(A,v)$ with $m(\le n)$ a positive integer, we often have $W$ swinging back and forth around the leading eigenspace $L_s(A)$ for some small integer $s$. Another problem with Krylov subspace technique is, in large scale computation it is impossible for us to keep all base vectors of $W$ when $A$ is not symmetric, even if the matrix $A$ might be sparse. Thus the projection of $x$ on subspace $W$ can not be obtained easily. Although in case $A$ is symmetric it is not necessary to keep all base vectors of  $W$ because of the three-term recurrence relations, in practical application we often encounter the problem of  so-called loss of orthogonality.

In view of (\ref{errorNorm}), our intention here is to find a vector $p$ so that the length $||p||$ of the projection vector of $x$ is as large as possible.  We start from an initial direction $p_0$ on which the projection of $x$ is known or easily available. A searching direction is then needed for the purpose of constructing a vector $p_1$ so that  $x$ has a larger projection on $p_1$ than that on $p_0$ in terms of vector length. For any searching direction $d$ we need to
 have the projection of $x$ on $d$ easily obtainable. An arbitrarily chosen direction vector $d$ can not be used since we
  don't have information about the inner product between $d$ and $x$. Fortunately we have a lots of vectors available from the system (\ref{eq:1}) since
$a_i x = b_i \, (i=1,2,\cdots, n).$ i.e., all row vectors $a_i$ in matrix $A$ can be used as our searching directions.

As a starting direction(it is not necessary though) it is thus a possible choice for  us to use $k$-th row vector of $A$ where the subscript $k$ is chosen so that
$$ {|b_k| \over ||a_k||} = \max_{1\le i \le n}{|b_i| \over ||a_i||}.$$
Yet a better starting direction seems to be $p_0=A^Tb$ (assuming $||A||_i = 1, i=1,2,\cdots, n $) since we have
$x^Tp_0 = p_0^Tx =b^TAx = b^Tb$ and hence a projection with larger length(i.e., ${b^Tb \over ||A^Tb||} > {|b_k| \over ||a_k||} $)
maybe available. The construction of next projection direction $p_1$ depends on a carefully chosen searching direction
 vector $d$ such that $|x^Td|/||d||$ is as large as possible.  There are many ways of determining a suitable
 searching direction $d$, however the following facts should be observed when we start the searching process.

Assume $x^Tv_i = b_i, \, (i=1,2$) with $b_1 \neq 0$ and $ ||v_i||=1
\, (i=1,2)$, we wish to find a real number $t$ such that the
function $f(t) $  defined by \be \label{def:f}
     f(t) = \frac{|x^Tv|}{||v||}
\ee
is maximized among all possible vectors in the form $v = v_1 + tv_2$.
It is easy to see from analysis that the answer to the above optimization problem lies on the following conclusion.

\begin{lemma}\label{lemma:3-1}
   Let $x^Tv_i = b_i, \, (i=1,2$) with $|b_1| \ge |b_2|$ and $ ||v_i||=1 \, (i=1,2)$,
and $ \alpha = v_1^Tv_2$. Let $ s = \frac{b_2 - \alpha b_1}{ b_1 - \alpha b_2 }$. Then
      \be \label{def:fs}
          f(s) \equiv \frac{|x^T(v_1+sv_2)|}{||v_1 + sv_2||}=\dstyle{\max_{t\in R}} \frac{|x^T(v_1 + tv_2  )|}{||v_1 + tv_2  ||}.
      \ee
Furthermore
\be
        f(s)\ge \max\{|b_1|,|b_2|\}
     \ee
\end{lemma}

{\it Proof.}  Let $$g(t) =\frac{x^T(v_1 + tv_2)}{||v_1 + tv_2||}.$$ We have
  $$  g(t) = \frac{b_1 + tb_2}{\sqrt{ 1 + 2\alpha t + t^2}}.$$ Thus
   $$
   \barr{ll}
   g'(t) &= \dstyle\frac{b_2  (1 + 2\alpha t + t^2) - (b_1 + t b_2)(\alpha + t)}{  (1 + 2\alpha t + t^2)^{3/2}} \\[4mm]
     & = \dstyle\frac{b_2 - \alpha b_1 - (b_1 - \alpha b_2)t } {(1 + 2\alpha t + t^2)^{3/2}  } \\
     &=\dstyle\frac{(b_1 - \alpha b_2)(s - t) }{(1 + 2\alpha t + t^2)^{3/2}  }
  \earr
  $$
  Let $g'(t)=0$ we have the solution as $t = \dstyle\frac{b_2 - \alpha b_1}{ b_1 - \alpha b_2 }\equiv s,$ i.e., $s$ is an
  extreme point for function $f(t)$. \\
  case 1. $ b_1> \alpha b_2,$ we have $g'(t)>0$ if ($t< s$) and $g'(t)<0 $ if $t>s$. That means $g(t)$ reaches the maximal
  value at $s$. Since $g(t) \rightarrow -b_2$ when $t \rightarrow - \infty$ and $g(t) \rightarrow b_2$ when $t \rightarrow +\infty$,
  we have $g(s)\ge g(t)>-b_2$ for all $t<s$ and $ b_2 < g(t)\le g(s)$ for all $t>s$, thus we have $f(t)=|g(t)|$ reaches its maximal
  value at $s$.\\
  case 2. $ b_1 < \alpha b_2,$, we have $g'(t)<0$ if ($t< s$) and $g'(t)>0 $ if $t>s$. That means $g(t)$ reaches the minimal value
  at $s$. Since $g(t) \rightarrow -b_2$ when $t \rightarrow - \infty$ and $g(t) \rightarrow b_2$ when $t \rightarrow +\infty$,
   we have $g(s)\le g(t)<-b_2$ for all $t<s$ and $ b_2 > g(t)\ge g(s)$ for all $t>s$, thus we have $f(t)=|g(t)|$ reaches its
    maximal value at $s$.\\
  Thus in both cases we have $f(s)> |b_2|$. Since $f(0) = | g(0)| = |b_1|$ and
  $f(s)$ is the maximal value of $f(t)$, thus we also have $f(s)> | b_1|$. See figure 1.

\begin{figure}[h] \label{fig:proj}
\begin{center}
\includegraphics[height=5cm,width=6cm]{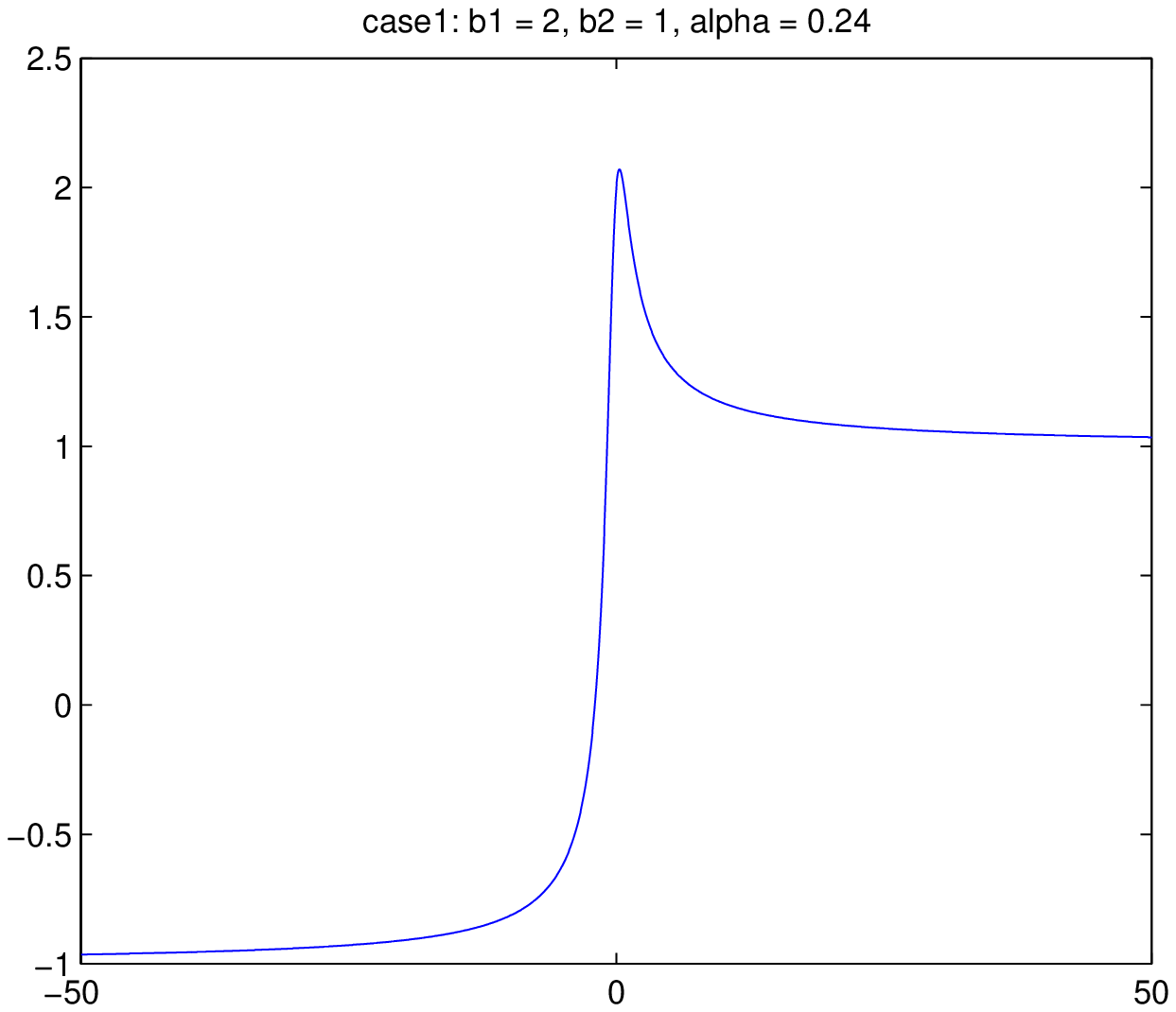}
\includegraphics[height=5cm,width=6cm]{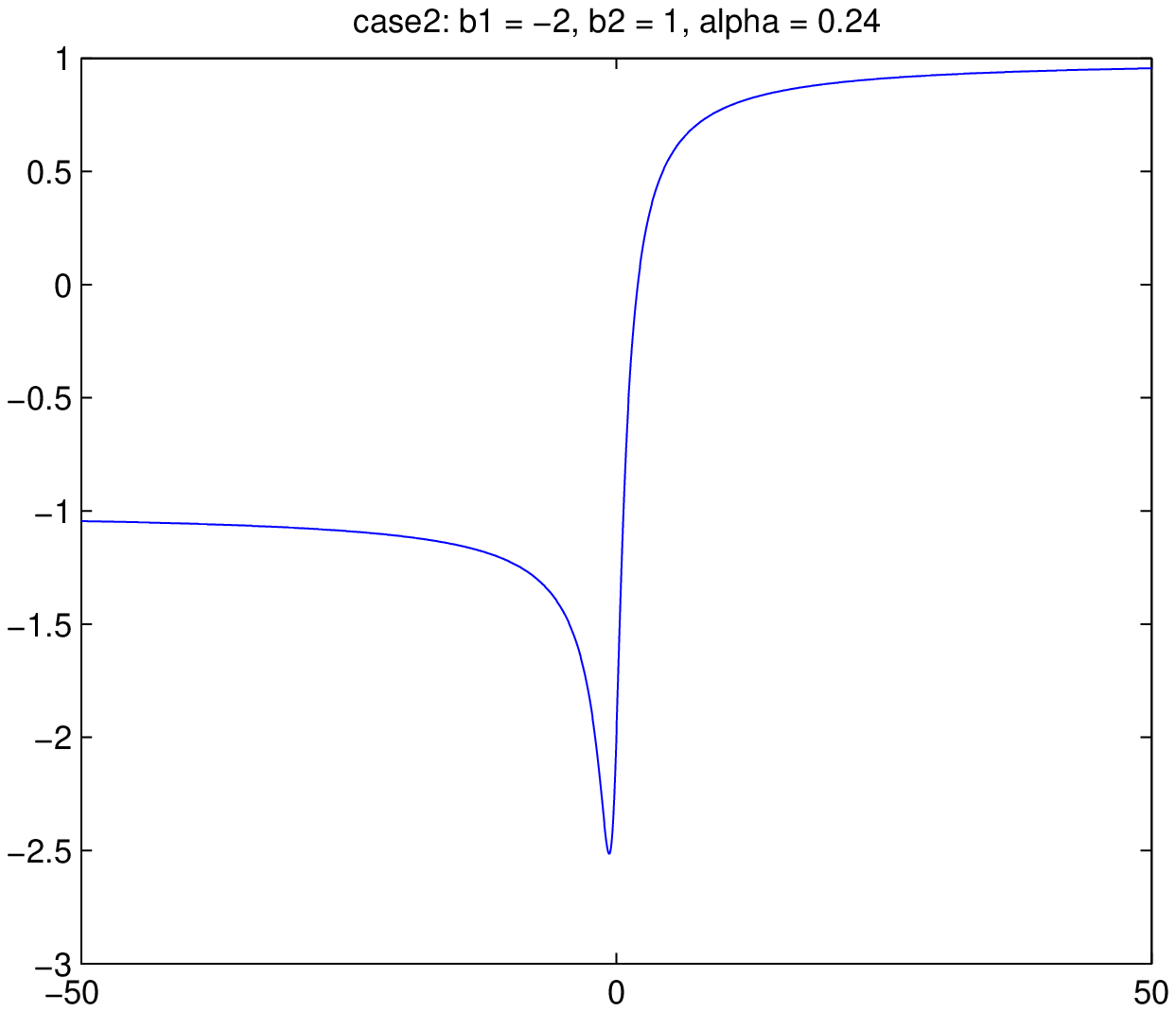}
\end{center}
\caption{Graph of g(t)}
\end{figure}

Remark: Assuming $b_1 \neq 0$, $f(s)$ can be rewritten as following by replacing $s$ as $s = \frac{b_2 - \alpha b_1}{b_1 - \alpha b_2}$.
\be \label{fs2}
    f(s)  = \dstyle\frac{|b_1 +  s b_2|}{\sqrt{ 1 + 2\alpha s + s^2}}
          = \dstyle\frac{\sqrt{ 1- 2\alpha r + r^2}}
             {\sqrt{ 1 - \alpha^2} }|b_1|
          = |b_1|\sqrt{1 + \dstyle\frac{(r- \alpha)^2}{ 1 - \alpha^2} }
\ee where $r = b_2/b_1$.

In view of (\ref{fs2}), $f(s) \rightarrow \infty$ when $\alpha \rightarrow 1$ (assuming $r$ independent of $\alpha$).
It is thus attempting for us to get the next projection $p_1$ of $x$ with much bigger length(and thus hopefully more closer to $x$) by carefully selecting suitable vector $d$ with $d^Tx=b_2$ and $ \alpha = p_0^T*d$ is as close as possible to $1$(i.e., the angle between $p_0$ and $d$ should be very small). However this seems to be very hard and thus we turn to an easier scheme to fulfill our task---we  will use subspaces on which projections of $x$ are easily available. For this purpose we now generalize our  conclusion in Lemma \ref{lemma:3-1} into  following statement.

\begin{lemma} \label{lemma:3-2}
   Let  $ x, v_i\in R^n\, (i=1,2,\cdots, m)$,
and $ W=span\{v_1,v_2,\cdots, v_m\}$. Let $ P_W(x) $ be the projection of $x$ onto subspace $W$. Then
     $$ \frac{v_*^Tx}{||v_*||} = \dstyle\max_{v \in W}  \frac{|x^Tv|}{||v||}$$
where $v_* = P_W(x)$.
\end{lemma}
Proof. Without loss of generality we can assume $ || x || = 1$.  By definition of angles between vectors  we have
   $$ f(v)= \frac{|x^Tv|}{||v||} = \frac{|x^Tv|}{||v||||x||} = |\cos<x,v>|$$
   where $<x,v>$ denotes the angle between vector $x$ and $v$. Obviously $f(v)$ reaches its maximum value if and only
   if $<x,v>$ is minimized, which is true only when $v$ lies on the projection of $x$ onto subspace $W$.

By using this result, one can always expect a searching direction $d$ on which vector $x$ has a bigger
projection length $||P_d(x)||$ than any vector in subspace $W = span\{v_1,v_2,\cdots, v_m\}$
with $x^Tv_i \, (i=1,2,\cdots, m)$ given.
Since we have $n$ vectors $a_i \,(i=1,2,3,\cdots,n)$ to form subspaces of $R^n$, this give us plenty of choices
 when it comes to construct subspaces. More importantly we can use parallel process to construct these
  subspaces and figure out projections of $x$ on each of them. Instead of using successive ``partial" projections
  which did not adequately make use of current system information, all these projections of
   $x$ can be used to construct a better approximation to the current system.

\subsection{The projection algorithms}
In this subsection we present some basic algorithms for solving
linear system of equations. We first introduce two algorithms for
calculating a projection vector $p$ of $x$ to the system
(\ref{eq:1}) based on current system data, i.e., the coefficient
matrix  $A$ and  right-hand side vector $b$, which is always the
unique vector in some subspace of $R^n$ on which solution vector $x$
having the maximum projection length.

In preparation, we begin with the division of all row vectors of $A$
into groups of vectors $\{G_i\}_1^k$, with each group $G_i$ contains
$m_i$ vectors, where $m_i \,(i=1,\cdots, k)$ are relatively small
integers satisfying $m_i<m, \, \forall  1 \le i \le k$. $m$ is a
suitable integer so that the QR factorization of matrix $T_i$ formed
by all vectors in group $G_i$ is applicable; in case of sparse
coefficient matrix, QS factorization process based on LGO
method \cite{pengDDM2009} can be used and thus $m$ can be
relatively large(say, up to $O(10^5)$). The right-hand side vector
$b$ is divided correspondingly into vectors $b_i\, (i=1,\cdots, k)$.

One thing needs to be mentioned here is that we assume two adjacent
groups $G_i$ and $G_{i+1}$ contain about half of their vectors in
common and any row vector in $A$ must lie in at least one of the
groups, we will refer this group $\{G_i\}$ as an overlapped division
of $A$. A non-overlapped division of $A$ means the intersection of
any two groups in the division is empty.

The first accumulated projection algorithm uses an overlapped division of $A$ and seeks the
projection vector $p_i$ of solution $x$ on the range of each  group
$G_i$, i.e., the subspace spanned by all vectors in $G_i$. All
projection vectors $\{p_i\}_1^k$ are then ``glued" together to form
a better projection vector of $x$, while the ``gluing" process is
nothing but another projection of $x$ over subspace $W =
span\{p_1,\cdots, p_k\}$. The details comes as follows.

\begin{alg}(AP version 1) \label{alg:AP1} Let $A \in R^{n\times n} $ be nonsingular, $b \in R^n$. The following
procedure produces a projection(vector) $p$ of solution $x$ to the
system $Ax = b$.
\begin{itemize}
\item{Step 1.}  Divide matrix $A$ into $k$ blocks:  $A  = [A_1^T, A_2^T, \cdots, A_k^T]^T$, divide $b$ correspondingly:
$b = (b_1^T, b_2^T, \cdots, b_k^T)^T$.

\item{Step 2.}  For each $i \in \{1,2,\cdots, k-1\}$, compute projection $p_i$ of $x$ in
$ran(A_i^T, A_{i+1}^T)$:  $p_i = G_i^T(G_iG_i^T)^{-1}B_i$ and
compute scalar $c_i = x^Tp_i$ as $c_i = B_i^T(G_iG_i^T)^{-1}B_i,$
where $G_i = [A_i^T, A_{i+1}^T]^T$ and $B_i = (b_i^T, b_{i+1}^T)^T$.

\item{Step 3.} Construct matrix $H$ as $H = [p_1, p_2, \cdots, p_{k-1}]$ and vector
$ c^T = (c_1, c_2, \cdots, c_{k-1})$.

\item Step 4. Form a projection $p$ of $x$ over $ran(H)$ and compute scalar \mbox{ $ \alpha $ ($= x^Tp$).}

\item Step 5. Output $p$ and $\alpha$.
\end{itemize}
\end{alg}

{\it Remark:} \begin{quote}
The projection process on each group of
row vectors can be handled independently and thus good for parallel
implementation.

There exists an important relation between $H$,$c$ and solution
vector $x$:  $$ H^Tx = c$$

In case the number of groups is too big so that a direct projection
over $ran(H)$ is not applicable, one can use a nested version of
this algorithm over $H$ to obtain the final projection vector.
\end{quote}

Algorithm \ref{alg:AP1} uses a sequence of projections on the a set of subspaces
determined by submatrices of $A$, these projections can be
obtained in parallel, which differs itself with those in Karcmarz's
idea. Furthermore, the blocks of matrices are overlapped with each
other. One can of course use different strategies when dividing the
matrix $A$ into submatrices and $b$ correspondingly.  It is easy to see
that direction vector $p$ satisfies
     $$\frac{x^Tp}{||p||} \ge \max_{i\in \{1, \cdots,n\} } \frac{a_ix}{||a_i||}$$
where $a_i$ denotes the $i$th row of matrix $A$. However $p$ may not
be the best option in general.

Another accumulated projection idea is to use a sequential projection process to get a
final projection vector $p$ of $x$. We begin with an initial
projection vector $p_0$ of $x$ and let it combine with all row
vectors in the first group $G_i$ to form a subspace $W_1$ of $R^n$,
and then find the projection vector $p_1$ of $x$ in $W_1$. $p_1$ is
then used to combine with all row vector in the next group $G_2$ to
form a subspace $W_2$ so that a projection vector $p_2$ of $x$ in
$W_2$ can be obtained. The above process is repeated until all
groups are handled so that the final projection vector $p_k$ is
available.   The following algorithm gives the details.

\begin{alg}\label{alg:AP2}(AP version 2) The following procedure
produces a projection vector of $x$ to system $Ax = b$.

\parbox{4mm}{~~}\parbox{11cm}{
\begin{itemize}
\item[step 1.]Divide matrix $A$ into $k$ blocks:
 $A  = [A_1^T. A_2^T, \cdots, A_k^T]^T$, divide $b$ correspondingly:  $b = (b_1^T, b_2^T, \cdots,
 b_k^T)^T$. Let $p_0 = \alpha A^T b$, $c = \alpha b^Tb$ where $\alpha =(b^Tb)/(b^TAA^Tb)$.
\item[step 2.] For $i =1$ to $k$  \\
\parbox{6mm}{~~}\parbox{10cm}{
\begin{itemize}
\item[step 2.1.] Construct matrix $W =[p_{i-1}, A_i^T] $ and vector $l =[c,b_i^T]^T$.
\item[step 2.2.] Compute the projection vector $p_i$ of $x$ onto subspace $ran(W)$  and the scalar $c(=x^Tp) $.
\item[step 2.3:] Go to next i.
\end{itemize}
}
\item[step 3:] Output $p_k$ and $c$.
\end{itemize}
}
\end{alg}

It is observed that Algorithm 2 is more effective than Algorithm 1 in terms of the
length of final projection vector $p_k$. By this reason, we use \mbox{Algorithm \ref{alg:AP2}} in our  numerical experiments(PAP and APAP algorithms).
It should be mentioned here that the AP algorithms depicts a successive projection
process over subspace $W_i = span\{p_i, v_1,v_2, \cdots, v_m\}$ ($i=1,\cdots, k$), where $v_1, v_2, \cdots, v_m$ denotes the row vectors of submatrix $A_i$, and $p_i=P_{i-1}x$ is the projection of $x$ over subspace $W_{i-1}$ with $p_0$ stands for the initial projection vector of $x$. Hence the whole AP process can be written in the matrix form as $p =  P_k x$ where $P_i$ ($i=1,\cdots ,k$) represents the projection matrix over subspace $W_i$. It is easy to see that $P_i$ depends on vector $x$.  As a matter of fact, $P_k$ has the form
\be \label{def:P_k}
   P_k = H(H^TH)^{-1}
\ee
where $H = [p_{k-1}, V_k]$ and $V_k = [v_1, v_2, \cdots, v_k]$, assuming $p_{k-1} \notin span{V}$.

As a straightforward application, Algorithm 1 and 2 can be used to
solve the linear system(\ref{eq:1}) as stated in the next algorithm.\\

\begin{alg}\label{alg:PAP}(Progressively Accumulated Projection Method--PAP). Let $A \in R^{n\times n}
$, $b \in R^n$. The following procedure produces an approximation
$y$ to the solution $x$ satisfying $Ax = b$.

\parbox{4mm}{~~}\parbox{10cm}{
\begin{itemize}
\item[step 1.] Initialize vector $y$ as zero vector.
\item[step 2.] While not converged

        \parbox{4mm}{~~}\parbox{10cm}{
        \begin{itemize}
        \item[step 2.1]   Use algorithm 1 or 2 to get a projection $p$ of $x$ to system $Ax = b$.

        \item[step 2.2]   Update $y$ as $y = y+p$.

        \item[step 2.3]   Update $b$ as $b = b - Ap$.

        \item[step 2.4]   Check convergence condition.
        \end{itemize}
        }
\item[step 3.] Output $p$
\end{itemize}
}
\end{alg}

PAP is based on the principles in section 2, hence the convergence (Theorem \ref{thm:convergence}) of this algorithm is straight-forward and its proof is thus omitted.
We need to mention here that unlike classical Krylov subspace methods, the AP-type methods proposed here can actually be used to solve any under-determined systems. Also we have to point out that each sweep in step 2 is a projection process with projection matrix $P_k$ varies. The following graph shows the comparison of
this algorithm at different iterative numbers, and Table \ref{table:PAPiter} gives the
needed iterations for a convergent solution under given tolerance, where
the coefficient matrix $A$ is chosen as $A = tridiag(-1,2,-1)$ with $A \in R^{100\times 100}$ and the block size is chosen as $20$ when applying algorithm 2 in this case.
\parbox{3cm}{~~}

\begin{figure}[h] \label{fig:proj}
\begin{center}
\includegraphics[height=5cm,width=10cm]{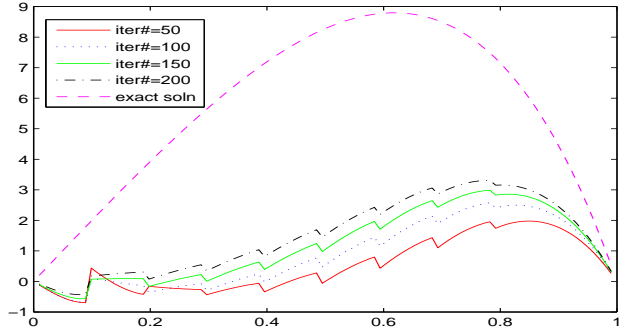}
\parbox{10cm}{\caption{Comparison of approx. solns at different \mbox{iteration} numbers}}
\end{center}
\end{figure}

\parbox{2cm}{~~}

\begin{table}\begin{center}
\caption{iteration numbers needed for convergence}
\label{table:PAPiter}
\begin{tabular}{|c|c|c|c|c|c|c|c|}
\hline  	 tolerance  & $10^{-1}$ &$ 10^{-2}$ &$ 10^{-3} $ & $10^{-4}$ & $10^{-5} $ &$10^{-6}$&  $10^{-7}$  \\ \hline
   iter\# & $3844$& $ 5534$& $ 7224$ &$8916$ &$10606$ & $12296$ & $13986$      \\ \hline
\end{tabular}
\end{center}
\end{table}

\section{Properties of AP process }

In this section we present some  analysis results for AP process described in Algorithm \ref{alg:AP2} .
\begin{lemma} \label{lemma:AP1}

 Assume that matrix $A \in R^{m\times n}$ ($m \le n$) has full row rank, $x \in R^n$ and $b \in R^m$ where $m\le
 n$ satisfying $Ax = b$.  Let $A$ be divided into $k$ submatrices by its rows:  $A = (A_1^T,A_2^T,\cdots, A_k^T)^T $
 with $A_i \in R^{m_i\times n}$, and $b$ is divided as  $b = (b_1^T,b_2^T,\cdots, b_k^T)^T$ correspondingly.   Let
 $\{p_i\}_1^k$ be the vector sequence produced by
 AP process(Algorithm 2).

\begin{enumerate}
      \item[(1)]    There holds for every $i=1,2,\cdots k$
           \be \label{eq:2}
            (x - p_i,p_s) = 0, \quad (s=i,i-1).
            \ee
      \item[(2)]
           Vector $ p_{i+1} - p_i$ ($i=0,1,\cdots, k-1$) is orthogonal to $p_i$, i.e.
           \be  \label{eq:3}
            (p_{i+1} -  p_i,p_i)=0
           \ee

      \item[(3)] There holds for $i=1,2,\cdots, k$
            \be  \label{eq:4}
              ||p_i||^2 + ||p_{i+1} - p_i||^2= ||p_{i+1} ||^2
            \ee

      \item[(4)]  For every $s (1 \le s \le k)$,  there holds
          \be  \label{eq:5}
            ||p_{s}||^2 = ||p_0||^2 + \sum_{i=1}^s {||p_i-p_{i-1}||^2}
          \ee
   \end{enumerate}

 \end{lemma}

{\bf\it Proof.}
  (1)  We first show that $(x-p_0,p_0) =0$. As a matter of fact, since $\alpha = b^Tb/(b^TAA^Tb)$,
  we have
     $$ (x - p_0,p_0) = (x - \alpha A_1^Tb_1, \alpha A^Tb) = \alpha x^TA_1^Tb - \alpha^2 b^TA_1A_1^Tb
            =  \alpha b^Tb - \alpha b^Tb = 0$$
  From the fact that $p_i$ is the projection of $x$ over subspace $ran(W_i)$ with
     $W_i = [p_{i-1}, A^T_i]$, for any $i ( 1\le i \le k)$ we must have
          \begin{center} $(x-p_i, p_i)=0 $ and $(x - p_i, p_{i-1})=0$
     \end{center}
      since  both $p_i$ and $p_{i-1} $ belong to $W_i$.

  (2) Note that from (\ref{eq:2}) we have
      $$ (p_{i+1} -p_i, p_i) =( ( x-p_i) - (x - p_{i+1}), p_i) =
             ( x-p_i, p_i) -(x - p_{i+1}, p_i)=0, $$
      which yields (\ref{eq:3}).

  (3) From (\ref{eq:3}) we have  $$
       \begin{array}{ll}
        || p_{i+1} - p_i||^2 &= ( p_{i+1}- p_i, p_{i+1} - p_i) \\
                              &=(p_{i+1}- p_i, p_{i+1} ) \\
                              & =  (p_{i+1}, p_{i+1}) - (p_i, p_{i+1}) \\
                              & =  (p_{i+1}, p_{i+1})  - (p_i, (p_{i+1} - p_i) + p_i)\\
                              & =    (p_{i+1}, p_{i+1}) - (p_i,p_i) \\
                              &  =  || p_{i+1} ||^2 - ||p_i||^2
       \end{array}
       $$
       from which (\ref{eq:4}) comes immediately.

  (4) Equation (\ref{eq:5}) follows from the recursive application of (\ref{eq:4}):
        $$
        \begin{array} {rl}
          ||p_s||^2 &= ||p_{s-1}||^2 + || p_s - p_{s-1}||^2 \\
                   & =||p_{s-2}||^2 +|| p_{s-1} - p_{s-2} ||^2 +   || p_s - p_{s-1}||^2 \\
                   & \cdots \\
                 & = ||p_0||^2 + ||p_1 - p_0 ||^2 + ||p_2 - p_1||^2 + \cdots +   || p_s - p_{s-1}||^2 .

        \end{array}
        $$ Proof is completed  \hfill$\Box$

Lemma \ref{lemma:AP1} actually tells the fact that the ``length"(norm) sequence $\{||p_i||\}_1^k$ of
projection vector $\{p_i\}_1^k$ actually forms a monotonically increasing sequence, and obviously $||x||$ is
actually
one of its upper bounds. In order to find out how fast this sequence is increasing, we need to figure out the
detailed information of each $||p_i|| \,(i=1,2,\cdots, k)$. The following conclusion answers this question.

\begin{lemma} \label{lemma:buildpi1}

 Assume the same assumption in Lemma \ref{lemma:AP1}. Let $\tilde{p}_i$ and $\tilde{x}_{i+1}$ be the projection vectors of $p_i$ and $x$ over
  subspace $ran(A_{i+1}^T)$ respectively. Then $p_{i+1}$ has the following expression
   \be \label{eq:6}
     p_{i+1} =  \alpha_i p_i + A^T_{i+1} u = \alpha_i (p_i - \tilde{p}_i) + \tilde{x}_{i+1}
   \ee
   where $u$ is
   \be\label{eq:7}
    u = \tilde{A}_{i+1} (b_{i+1} - \alpha_i A_{i+1} p_i)
    \ee
   and
   \be \label{eq:8}
    \alpha_i = \frac{x^Tp_i - (A_{i+1}p_i)^T \tilde{A}_{i+1}  b_{i+1}    }
   { p_i^Tp_i - p_i^TA_{i+1}^T\tilde{A}_{i+1} A_{i+1} p_i }
   = \frac{x^Tp_i - p_i^T*\tilde{x}_{i+1}}{p_i^Tp_i - \tilde{p}^T_i\tilde{p}_i}
   \ee
   and
   $$ \tilde{A}_{i+1} = (A_{i+1} A_{i+1}^T)^{-1}.$$
   Furthermore
    \be\label{eq:norm_pk} ||p_{i+1}||^2 = \alpha_i ||p_i||^2  + b_{i+1}^T \tilde{A}_{i+1} b_{i+1}-
    \alpha_i^2(A_{i+1}p_i)^T\tilde{A}_{i+1} (A_{i+1}p_i)
    \ee

 \end{lemma}
{\bf \it Proof.} \\
 It is valid to express $p_{i+1}$ in the form like the first equation of (\ref{eq:6}) for some $u \in R^{m_{i+1}}$  since $ p_{i+1} \in W_{i+1} = ran([p_i, A_{i+1}^T]) $, where $m_i$ is the number
 of rows in submatrix $A_i$.

Since $p_{i+1}$ is the projection of $x$ over subspace $W_{i+1}$, we have
    $$ A_{i+1}(x-p_{i+1})=0$$
which leads to
    $$ b_{i+1} - \alpha_i A_{i+1} p_i - A_{i+1}A_{i+1}^T u = 0$$
from which comes (\ref{eq:7}).
 Note that
 $\tilde{p}_i $ and $\tilde{x}_{i+1}$ are projections of $p_i$ and $x$ over $ran(A_{i+1}^T)$,
 we have
 \be\label{eq:xx}
  \tilde{p}_i = A_{i+1}^T \tilde{A}_{i+1} A_{i+1}p_i\quad\mbox{and}\quad\tilde{x}_{i+1} = A_{i+1}^T \tilde{A}_{i+1} b_{i+1}
  \ee
Plug (\ref{eq:7}) and (\ref{eq:xx}) back into the first equation of (\ref{eq:6}) gives the second equation
of (\ref{eq:6}).

 Similarly, by $ p_i^T(x-p_{i+1}) =0$ we have
     $$ x^Tp_i - \alpha_i p_i^Tp_i + u^TA_{i+1}p_i =0,$$
 replacing $u$ by  (\ref{eq:7}) yields the first equation of (\ref{eq:8}).
 Again because $\tilde{p}_i$ is the projection of $p_i$,
 $$
   (p_i-\tilde{p}_i, \tilde{p}_i) = 0,
 $$
 this means
 \be \label{eq:yy}
   p_i^T \tilde{p}_i = \tilde{p}_i^T \tilde{p}_i.
 \ee Plug (\ref{eq:xx}) and (\ref{eq:yy}) into the first equation of (\ref{eq:8}) gives the second equation.

 Finally from (\ref{eq:6}) we have
 \be     \label{eq:9}
     \barr{rl}
        ||p_{i+1} ||^2 & = (\alpha_i p_i + A_{i+1}^T u)^T(\alpha_i p_i    + A_{i+1}^T u)  \\
                       & = \alpha_i^2 p_i^Tp_i + 2 \alpha_i p_i^TA_{i+1}^T u+ u^TA_{i+1}A^T_{i+1}u.   \\

     \earr
 \ee
 Since
 \be \label{eq:10}
 \barr{rl}
           2 \alpha_i p_i^TA_{i+1}^T u &= 2 \alpha_ip_i^TA_{i+1}^T\tilde{A}_{i+1}(b_{i+1} - \alpha_i A_{i+1}^Tp_i)\\
                                     & = 2 \alpha_ip_i^TA_{i+1}^T\tilde{A}_{i+1}b_{i+1} - 2\alpha_i^2
                                     (A_{i+1}p_i)^T\tilde{A}_{i+1}(A_{i+1}p_i)
   \earr
 \ee
 and
\be        \label{eq:11}
  \barr{rl}
            u^TA_{i+1}A^T_{i+1}u &= ( b_{i+1} - \alpha_i A_{i+1}p_i)^T\tilde{A}_{i+1}( b_{i+1} - \alpha_i
            A_{i+1}p_i)\\
                               & = b_{i+1}^T\tilde{A}_{i+1} b_{i+1} - 2\alpha_i b_{i+1}^T \tilde{A}_{i+1}
                               A_{i+1}p_i \\
                               & ~~+\alpha_i^2 (A_{i+1}p_i)^T\tilde{A}_{i+1}(A_{i+1}p_i), \\
     \earr
 \ee
equation (\ref{eq:norm_pk}) comes from (\ref{eq:9})  (\ref{eq:10}) (\ref{eq:11}) combined. \hfill$\Box$

Lemma \ref{lemma:buildpi1} describes one way of constructing $p_{i+1}$, and detailed information about
$p_{i+1}$
is revealed by (\ref{eq:norm_pk}). However a more direct approach can be used to evaluate the difference of the
norms between two  consecutive projections $ p_{i+1}$ and $p_i$. These can be shown in the following
conclusion.

\begin{lemma} \label{lemma:buildpi2}
 Assume the same assumption in Lemma \ref{lemma:AP1}. Let $I$ be the identity matrix in $R^n$. Then $p_{i+1}$
 has the following expression
\be\label{eq:piform2}
    p_{i+1} = p_i + \bar{A}_{i+1}^T v
\ee
and
\be \label{eq:pinorm}
 ||p_{i+1}||^2 - ||p_i||^2   = (b_{i+1} - (p_i^Tx)d)^T(\bar{A}_{i+1}  \bar{A}_{i+1}^T)^{-1}  (b_{i+1} - (p_i^Tx)d)
\ee
where
   $\bar{A}_{i+1}$ is a rank-one modification of submatrix $A_{i+1}$ as
   \be\label{def:Ai_bar}
      \bar{A}_{i+1}=  A_{i+1} - d p_i^T = A_{i+1}(I - u_iu_i^T)
   \ee
   with $u_i = p_i/||p_i||$,  $d \in R^{m_{i+1}}$ a vector taken as
   $ d =  A_{i+1}p_{i+1} /  ||p_i||^2 $
   and $v$ is defined as
   $$
      v = (\bar{A}_{i+1}\bar{A}_{i+1}^T)^{-1} (b_{i+1} - (x^Tp_i)d )
   $$
   assuming the related inverse exists.
\end{lemma}
\emph{Proof.} \\
Since $p_{i+1}$ is the projection of $x$ over subspace $W_{i+1} = ran([p_i, A_{i+1}^T])$,
it can be constructed as follows.

First we modify row vectors in $A_{i+1}$ so that they are orthogonal to vector $p_i$, this
can be depicted as a rank-one modification to $A_{i+1}$ as
   $$ \bar{A}_{i+1} = A_{i+1} - d p_i^T,$$
where $d$ can be obtained from the fact that
  $$ \bar{A}_{i+1} p_i = 0$$
which leads to
   $$ A_{i+1} p_i - d p_i^Tp_i =0,$$
hence
  $$ d = A_{i+1}p_i/(p_i^Tp_i), $$ and
  $$\bar{A}_{i+1} = A_{i+1} - dp_i^T = A_{i+1} - A_{i+1}p_ip_i^T/(p_i^Tp_i) = A_{i+1}(I - u_iu_i^T),$$
  where $u_i = p_i/||p_i||$.

Next we calculate the projection vector $\tilde{p}_{i+1}$  of $x$ over $ran(\bar{A}_{i+1})$ as
$$ \tilde{p}_{i+1} = \bar{A}_{i+1}^Tv,$$
where $v$ can be derived from the fact that
$$ \bar{A}_{i+1} (x - \tilde{p}_{i+1} ) = 0,
$$
which  leads to
$$
   v = (\bar{A}_{i+1}  \bar{A}_{i+1}^T)^{-1} (b_{i+1}- (p_i^Tx)d )
$$
assuming    $(\bar{A}_{i+1}  \bar{A}_{i+1}^T)^{-1}$ exists.

Since $\tilde{p}_{i+1} = \bar{A}_{i+1} v$ is the projection of $x$ over $ran (\bar{A}_{i+1})$ and
$\bar{A}_{i+1} p_i = 0$, we must
have $(p_i, \tilde{p}_{i+1}) =0 $. Therefore
$$
\barr{rl}
   ||p_{i+1}||^2 - ||p_i||^2 &= ||\tilde{p}_{i+1} ||^2             \\
                             & = v^T\bar{A}_{i+1} \bar{A}_{i+1}^Tv \\
                             & = (b_{i+1} - (p_i^Tx)d)^T(\bar{A}_{i+1}  \bar{A}_{i+1}^T)^{-1}  (b_{i+1} -
                             (p_i^Tx)d).
\earr
$$
noting that matrix $(\bar{A}_{i+1}  \bar{A}_{i+1}^T)^{-1}$ is symmetric(actually positive definite
symmetric).\hfill$\Box$.

\emph{Remark}: It can be shown that the length difference between $p_{i+1}$ and $p_i$ can also be written as
\be \label{diff:len}
    ||p_{i+1}||^2 - ||p_i||^2 = \tilde{x}^T G \tilde{x}
\ee
where $ G =(\bar{A}_{i+1}\bar{A}_{i+1}^T)^{-1}$ and $\tilde{x} = \bar{x} - (x^Tu)u$, where $\bar{x}$ denotes the
projection of $x$ on $ran(A^T)$ and $(x^Tu)u$ is the projection of $x$ (as well as $\bar{x}$) on the direction of
$u= p_i/||p_i||$.

 Note that in the above lemma, we need to assume the existence of each matrix  $(\bar{A}_{i+1}
 \bar{A}_{i+1}^T)^{-1}$. The following conclusion gives the sufficient and necessary conditions for these to hold
 true.

 \begin{lemma}\label{lemma:3-5}
    Let $A \in R^{m\times n} (m \le n)$ and $rank(A) = m$, $u \in R^n$ be a unit vector
    in $R^n$.  Let $\bar{A} = A(I - uu^T)$ and $G = \bar{A}\bar{A}^T$, where $I$ denote the identity
    matrix in $R^n$. Then $G$ is nonsingular if and only if $u \notin ran(A^T)$.
 \end{lemma}
 \emph{Proof}.  \\
 Note that $G=\bar{A}\bar{A}^T$ is invertible if and only if $ \bar{A}$ is of full row rank.

       (Necessity) Assume $G$ is invertible, we need to show that $u \notin ran(A^T)$. If this is not the case,
       i.e., $u \in ran(A^T)$, then there is a $v \in R^n$ ( $v \neq 0$) such that $ u = A^Tv$.
       Thus
       $$
         \bar{A}^T v = (A(I- uu^T))^Tv = (A - Auu^T)^Tu = A^Tv - uu^TA^Tv = u - u(u^Tu) = 0
       $$
       since $||u||=1$.  This means $\bar{A}$ is not of full rank, hence $G$ is singular, a contradiction
       with our assumption.

       (Sufficiency). Suppose $u\notin ran(A^T)$, we need to show that $G$ is invertible. As a matter of fact, if
       $G$  is not invertible, then $\bar{A}$ is not of full-row rank. Therefore there exists a nonzero vector
       $v\in R^m$
       such that $\bar{A}^Tv =0$. That means
        $$ 0= (A(I - uu^T))^Tv = A^Tv - uu^TA^Tv = A^Tv - \alpha u $$
        where $\alpha = u^T(A^Tv)$ is a scalar. It is easy to see from here that $\alpha \neq 0$, otherwise we
        would have
        $A^Tv=0$ which means $A$ is not of full row rank. Hence $u = A^Tv/\alpha$, i.e., $u\in ran(A^T)$, this is
        contradictory    with the assumption.  \hfill$\Box$

\begin{lemma} \label{lemma:AP4}
 Assume the same assumption in Lemma \ref{lemma:AP1}. Vector sequence $p_0$,$p_1$,$\cdots$, $p_k$ are produced
 in
 one AP process, then
 \be \label{eq:pinorm_relation}
      ||p_i|| \le||p_{i+1}|| \quad ( i=0,1,2,\cdots, k)
 \ee
 and the equal sign holds  if and only if
   $$ A_{i+1} p_i = b_{i+1}$$
 \end{lemma}
\emph{ Proof. } Inequality (\ref{eq:pinorm_relation}) comes from (\ref{eq:4})  directly.
  We now prove the necessary condition for $||p_{i+1}|| = || p_i||$.

 \noindent (Necessity)
  Note that if $||p_{i+1}|| = ||p_i||$ holds , by (\ref{eq:4}) we must have $p_{i+1} = p_i$.
  Also from (\ref{eq:6}) we know that
    $$ p_{i+1} = \alpha_i p_i + A_{i+1}^T u, $$
  thus
  \be \label{eq:tt}
     A_{i+1}^T u = p_{i+1} -  \alpha_i p_i  = (1 - \alpha_i) p_i.
  \ee
  Multiplying both sides of (\ref{eq:tt}) by $A_{i+1}$ we have
  \be \label{eq:tt1}
     A_{i+1} A_{i+1}^Tu = (1- \alpha_i) Ap_i .
  \ee
  Note that from (\ref{eq:7}) we have
  \be \label{eq:tt2}
     A_{i+1} A_{i+1}^Tu =  b_{i+1} - \alpha_i Ap_i.
   \ee
  Combining (\ref{eq:tt1}) and (\ref{eq:tt2}) yields $$ Ap_i = b_{i+1}.$$

  \noindent(Sufficiency)Now we prove  $p_{i+1} = p_i$ under the assumption $ Ap_i = b_{i+1}.$

  As a matter of fact, in view of (\ref{eq:7}) and (\ref{eq:8}) we only need to show that
   $$ \alpha_i  = 1 $$
   in this case.

   Since  $ (x - p_i,p_i) = 0,$ we have
    \be\label{eq:tm1}
      x^Tp_i = p_i^Tp_i .
    \ee
    By using $ A_{i+1}p_i = b_{i+1}$ we obtain
    \be\label{eq:tm2}
       p_i^TA_{i+1}^T \tilde{A}_{i+1} A_{i+1}p_i = b_{i+1} \tilde{A}_{i+1} b_{i+1}
    \ee
    Hence from (\ref{eq:8}) we have
    $$
      \alpha_i = \frac{x^Tp_i - (A_{i+1}p_i)^T \tilde{A}_{i+1} b_{i+1} }
                      { p_i^Tp_i -  p_i^TA_{i+1}^T \tilde{A}_{i+1} A_{i+1}p_i}
               = \frac{p_i^Tp_i - b_{i+1}^T \tilde{A}_{i+1} b_{i+1} }
                      { p_i^Tp_i - b_{i+1}^T \tilde{A}_{i+1} b_{i+1} }
                =1
                $$
    This completes the proof of the sufficient condition. \hfill $\Box$

\section{An Accelerative Scheme}

We have observed from the preceding sections that the convergence
speed of the simple iterative algorithm may not be very satisfactory
in general. In this section we are to design some accelerative
approach for the PAP algorithm.

If we check the PAP procedure (Algorithm \ref{alg:PAP}) carefully and let $p_i, c_i$
denote the output from each call to algorithm \ref{alg:AP1} or \ref{alg:AP2}, the sum of
$p_i$ is used in Algorithm 3 as an approximation $x_k$, i.e., $x_k =
 p_1 + p_2 + \cdots + p_k$ when the $x_k$ satisfies some convergence
conditions, it is used as the final output approximation. The
following facts are obvious.

\begin{thm} \label{thm:convergence}
  Let $x_k$ be defined as above, $e_k = x - x_k$, then
$$\lim_{k\rightarrow \infty} x_k = x$$ i.e.,
$$\lim_{k\rightarrow \infty} e_k
= \lim_{k \rightarrow \infty} x - x_k = 0.$$ and
$\{||e_k||\}_1^\infty$ is strictly decreasing.
\end{thm}

We aim to use a combination of some selected approximate solution
from the sequence $\{x_k\}_1^M$ ( $ M \!< \!\infty$) to obtain an
``optimal" approximation, named as $y_1$, and then use it to update
the original system. This process is repeated to get next optimal
approximation solution $y_2$, and so on and so forth.

In Krylov subspace methods, this can be accomplished by using
coefficients of some specially chosen polynomials, as explicitly
done in Chebyshev semi-iterative method for accelerating stationary
methods or implicitly done in GMRES, FOM, etc. The drawbacks of
these techniques are: the resulted combination still falls into a
Krylov subspace with the same fixed generator matrix and starting
vector, which leads to their ultimate inefficiency in solving large
scale problems and thus have to resort to some preconditioning
techniques. Furthermore this treatment can not always guarantee a
convergent scheme.

Our strategy here is to pick up a subsequence of $\{x_k\}_1^M$, say,
starting from an initial $x_{k_0}$ in $\{x_k\}_1^M$ and then pick
another vector $x_{k_i}$ after every $t$ iterations to form a
sequence $\{x_{k_i}\}_{i=0}^m$ with $m$ a small integer. For
example, in the sequence $\{x_k\}_1^M$ we pick up $x_{10}, x_{20},
\cdots, x_{m*10}$ as the subsequence, renamed as $\{v_k\}_1^{m}$ and
then we try to find the projection of $x$ on the subspace $W =
span\{v_1, v_2, \cdots, v_m\}$.  To reach this goal we have to
resolve two key problems: first one should be able to obtain the
inner product between $x$ and each $x_k$; secondly one should be able to
selectively store the wanted subsequence $\{v_k\}_1^m$ from the
sequence $\{x_k\}_1^M$ and discard the unwanted vectors in the
sequence without affecting the calculation of inner product between
$x$ and $x_k$.

The following conclusion helps to resolve these questions.

\begin{thm} Let $x$ be the solution to system (1), $r_0 = b$, $p_k$ be the
approximation to $e_{k-1}$ in system $A e_{k-1} = r_k$ by Algorithm
1 or 2, $c_k = e_{k-1}^T p_k$ with $e_{-1} = x$, $r_k = r_{k-1}
- A p_{k-1}$, $e_k = e_{k-1} - p_k$, $x_k = \sum^k_{i=0} p_i$.  Then
\begin{align}
  x &= x_k + e_k, \label{thm1:1}\\
 x^T p_k &= c_k + x^T_{k-1} p_k \label{thm1:2}\\
 x^T x_k &= \sum^k_{i=0} c_i + \sum^k_{i=0} x_{i-1} p_i
 \label{thm1:3}
\end{align} for $k = 0, 1,2, \cdots$.
\end{thm}
{\it Proof.} Since $e_i = e_{i-1} - p_i$, we have
   $$e_{i-1} = e_{i} + p_i.$$
Hence
    $$ \sum_{i=0}^k e_{i-1} = \sum_{i=0}^k e_{i}+ \sum_{i=0}^k p_i,$$
which can be rewritten as
    $$\sum_{i=0}^k e_{i-1}- \sum_{i=0}^k e_{i} = \sum_{i=0}^k p_i.$$
Note that $e_{-1} = x$, we have
$$  x_k = \sum_{i=0}^k p_i =\sum_{i=0}^k e_{i-1}- \sum_{i=0}^k e_i = e_{-1} - e_k= x - e_k. $$
Thus
$$ x = e_k + p_k,$$
which gives (\ref{thm1:1}).

Note that
$$ 
   x^T p_k = (x_{k-1} + e_{k-1})^T p_k \\
            = x^T_{k-1} p_k + e^T_{k-1} p_k\\
            = c_k + x_{k-1}^T p_k
$$
which yields (\ref{thm1:2}).
Finally by (\ref{thm1:2}) we have
$$\barr{l}
    x^T x_k = x^T \sum_0^k p_i = \sum_0^k x^T p_i = \sum_0^k(c_i +
    x_{i-1}^T p_i) = \sum_0^k c_i + \sum_0^k x_{i-1}^T p_i
    \earr
$$
which gives (\ref{thm1:3}).  \hfill $\Box$

Expression (\ref{thm1:3}) suggests us that the inner produce between $x$
and $x_k$ for any $k$ only depends on two real number sequences
$\{c_i\}_0^k$ and $\{\tau_i\}_0^k$ with $\tau_i = x_{i-1}^Tp_i$ and
$\tau_i$ only depends on the last approximation $p_i$ to error
vector $e_i$ and accumulated approximation $x_{i-1}$. Hence it is
possible for us to design an algorithm which only needs to store a
constantly updating vector $x_k$ and save two number sequences
$\{c_i\}$ and $\{\tau_i\}$ during the iteration process. It is thus
viable for us to selectively store the wanted approximation $x_k$
and discard those unwanted ones in the approximation sequence
$\{x_i\}_1^M$.

The following algorithm makes use of these benefits and constantly
seeks a projection vector on the subspace formed by the selected
approximation vectors.

\begin{alg}\label{alg:APAP}(Accelerated Progressively Accumulated Projection--APAP) Let $x$ be the
solution to system (\ref{eq:1}), $\Delta \subset \{1,2,\cdots, M\}$
be a predetermined index set. The following procedure produces an
approximation to the solution $x$ in $Ax =b$.\\
\parbox{5mm}{~~}\parbox{11cm}{
    Step 1: (Initializing) Set $y\in R^n$ as zero vector\\
    Step 2: Do while not convergent\\
       \parbox{5mm}{~~} \parbox{10cm}{
        step 2.1  Set $r = b(\in R^n)$, $l_i = 0(\in R)$, $x_i = 0 \,(\in
        R^n)$\\
        step 2.2 For $i =1$ to $M$\\
            \parbox{5mm}{~~} \parbox{12cm}{
             step 2.2.1 Call Algorithm 1 or 2 to get projection vector $p_i$
             to $e_{i-1}$ satisfying $A e_{i-1} = r$ and $c_i = e_{i-1}^T
             p_i$.\\
             step 2.2.2 Set $\tau_i = x_i^T p_i$, $ l_i = l_i + \tau_i + c_i$\\
             step 2.2.3 Set $x_i = x_i + p_i$\\
             step 2.2.4 Update $r$ as $ r = r - Ax_i$.\\
             step 2.2.5 Store $x_i$ into matrix $H$ as a row vector and $l_i$ into vector $L$ if $i \in \Delta$.
             }\\
        step 2.3  Calculate projection vector $v$ of $\bar{x} (= A^{-1}b)$ on $ran\{H\}$ \\[-3mm]
            $$ v = H(H^TH)^{-1}L$$\\[-6mm]
        step 2.4  $y = y+v$\\
        step 2.5  $b = b - Av$
        }\\
    end
}
\end{alg}

{\it Remark:} \begin{quote}
In both PAP and APAP methods one has to repeatedly call basic AP methods(version 1 or version 2) to get the projections. Hence in actual implementation of these two algorithms it is necessary to rewrite the original system into its equivalent forms. In case the division of $A$ and $b$ is non-overlapped, each subsystem $A_i x = b_i$ corresponding to the division can be rewritten as $ Q_i^T x = \tilde{b}_i$ where $Q_iR_i=A_i^T$ forms the QR factorization(in case $A$ is dense) or QS factorization(in case $A$ is sparse) of $A_i$, while $  \tilde{b}_i= (R_i^T)^{-1} b_i$; in case of an overlapped division, one can use some extra sequence of submatrix-vector pairs to get the projections easily. By these rearrangement it is thus very efficient for us to get the projections of any vector on each subspaces. Note that the orthogonalization of each submatrix is needed only once.
\end{quote}

It turns out that the accelerating effect of this algorithm is remarkable by comparing Table \ref{table:PAPiter} and Table \ref{table:APAP1} where the block sizes in both tests are  exactly the same and the predetermined index set is selected as $\Delta =\{10,20,30,40,50,60\}$. For instance, to reach the same level $O(10^{-7})$ of relative residual error, PAP method(Algorithm \ref{alg:PAP}) needs almost 14000 iterations while \mbox{APAP} method(Algorithm \ref{alg:APAP}) needs only $2\times 60=120$ iterations, an amazingly improved convergence speed!

\parbox{2cm}{~~}

\begin{table}\begin{center}
\caption{iteration numbers needed for convergence using APAP}
\label{table:APAP1}
\begin{tabular}{|c|c|c|c|}
\hline  	 tolerance range  & $10^{-1}$ -- $ 10^{-7} $ & $10^{-8}$--$10^{-13} $ &$10^{-14}$--$10^{-19}$  \\ \hline
   outter iter\# & $2$ &  $3$ &$4$      \\ \hline
\end{tabular}\end{center}
\end{table}

\section{Numerical Experiments}

In this section we will show some applications of the aforementioned
\mbox{APAP} method. APAP is used to compare with block Jacobi method and GMRES since both are currently benchmark iterative methods in the category of extended Krylov subspace methods: the former is stationary and the later is non-stationary.

In the first example, we chose coefficient matrix $A\in R^{400 \times 400}$ as the following tridiagonal matrix
$$  A = tridiag(-1,2,-1)$$
and the solution vector $x$ is taken as  the values of function $u(t)=t(1-t)e^{3+t}$ at grid points $t = i\cdot h$ ($ i=1,2, \cdots, 400$) and $h = 1/401$. The results are listed in Table \ref{table:JacobiViaAPAP}.

\begin{table}[h]\begin{center}
\caption{Comparison between APAP and block Jacobi}
\label{table:JacobiViaAPAP}
\begin{tabular}{|c|c|c|c|c|c|c|}
\hline block size  & \multicolumn{2}{|c|}{cpu time(s)}  & \multicolumn{2}{|c|}{iter \#} &
\multicolumn{2}{|c|}{rel. residual} \\ \hline
  & blk Jacobi & apap & blk Jacobi & apap & blk Jacobi & apap   \\\hline
   30  &     34.9 &      2.70    &    11015   &       540&  4.03e-5&  1.59e-9 \\ \hline
   35   &    26.8  &     2.31    &     9528   &       440 & 3.73e-5 & 5.52e-11 \\ \hline
   40  &     20.2  &     1.60    &     8406   &       330 & 3.49e-5 & 1.38e-10 \\ \hline
   45  &     18.2  &     1.08    &     7533   &       220 & 3.29e-5 & 6.67e-10 \\ \hline
   50  &     16.4  &      1.36    &     6827   &       320 & 3.12e-5 & 4.27e-11  \\ \hline
\end{tabular}
\end{center}
\end{table}
We can see from this table that APAP exhibits much better performance than block Jacobi method does in terms of precision measured by the relative residuals, in the mean time APAP used much less cpu time and iteration numbers either.

As the second example, we use APAP to solve the Poisson problem defined on the unit square
$[0,1] \times[0,1]$.  The discretization scheme is the FDM five-point stencil, the resulted
coefficient matrix $A$ is a symmetric block diagonal matrix $A  \in R^{n\times n}$ ($n=2000$)
and the exact solution $u\in R^n$ is taken as grid values of function $u = x(1 -x)y(1-y)e^{3+x^2+y^2}$
at grid nodes $\{(i\cdot hx, j\cdot hy)\}$ with $i=1,\cdots, 50, j = 1,\cdots, 40$ and $hx = 1/51, hy = 1/41$.
The following table shows the iterations need for convergence with tolerance set as $10^{-5}$ as well as the comparison between the relative errors obtained by these two methods. Note that here the matrix $A$ has a relatively small condition number $cond(A) = 867$ and the block size is  determined as $ \sqrt{m n}$ so that an AP iteration needs approximately the same amount of storage as those of GMRES, where $n$ is the size of the system and $m$ is the predetermined restart number for $GMRES(m)$, note that the actual iteration number of $GMRES(m)$ is $out*in$ with $out$ as the specified maximum iteration for GMRES and $in$ as the restart number Matlab actually used in its running, while the actual iteration number is also counted as $out*in$. It seems that GMRES outperforms APAP in terms of time and iteration numbers in this case. We also need to mention here since the coefficient matrix $A$ is SPD, thus CG can be used to solve this problem and we recorded that CG outperforms both GMRES and APAP in this example in both cpu time and accuracy.

\begin{table}[h]\begin{center}
\caption{Comparison between APAP and GMRES}
\label{table:iters5}
\begin{tabular}{|c|c|c|c|c|c|c|c|c|}
\hline
\multicolumn{2}{|c|}{settings} & \multicolumn{2}{|c|}{iter. \#}  &  \multicolumn{2}{|c|}{time(in s)} &
\multicolumn{2}{|c|}{rel. error} \\ \hline
apap &gmres & apap & gmres & apap & gmres & apap & gmres \\\hline
 blk\_size& restart  &(out,in)  &(out,in)  &  &   &   &  \\ \hline
90&4 & (12,50) & (234,1)  &   13.3 & 8.4  &  7.6e-5  &  1.6e-4   \\ \hline
110&6 & (8,50) &(106,3)  &   9&   2.7 & 2.5e-5  &  1.6e-4\\ \hline
127&8 & (4,50) &(61,3 )   &   4.6& 3.9  & 7.5e-5 &   1.6e-4\\ \hline
142&10 & (4,50) & (40,4 )  &   4.6& 3.2  & 6.5e-5  &   1.6e-4\\ \hline
155&12 &(3,50) &(28,11)  &  4.1 & 3.1  &  7.6e-5  &  1.5e-5\\ \hline
168&14 &(3,50) &(22,2 )  &  4.1& 2.9   & 3.7.e-5  &  1.5e-5\\ \hline
179&16  &(3,50) &(17,11 ) &  4.1& 2.5  &5.1e-5 &   1.5e-5\\ \hline
190&18  &(2,50) &(14,13) & 3.0 & 3.4 &  5.2e-6  & 1.4e-5\\ \hline
\end{tabular}
\end{center}
\end{table}

The third test is on a system with asymmetric coefficient matrix $$A=tridiag(-1,2,-1.05)$$ having condition
numbers varying from $18944$ to $ 4.0902*10^{51}$ with $n$ varying from $100$ to $4600$, the following table shows the comparison between APAP and GMRES applied on the same systems with exact solution as $u = 2\sin(\pi x)e^{3+x}$. Note that the iteration number for APAP and GMRES are the total iteration numbers computed as inner loop multiplied by outer loop numbers. The restart number ($m$) for GMRES is fixed at $8$. It is interesting to see that the relative error of APAP is much better than that of GMRES, different than that in the second example.
\begin{table}[h]\begin{center}
\caption{Comparison between APAP and GMRES}
\label{table:iters5}
\begin{tabular}{|c|c|c|c|c|c|c|c|c|c|c|}
\hline
\multicolumn{1}{|c|}{apap} & \multicolumn{2}{|c|}{iter. \#}  &  \multicolumn{2}{|c|}{time(in s)} &
\multicolumn{2}{|c|}{rel. error} &\multicolumn{2}{|c|}{rel. residual }\\ \hline
size & gmres& apap & gmres &apap & gmres &apap & gmres &apap\\\hline
 29  &      8.0e+3    &  729    &    1.4 &     0.14 &    3.14e-2  & 6.71e-8  & 4.84e-4 & 1.37e-6\\ \hline
 70  &      4.8e+4   &  2016    &   10.3  &    0.95 &  8.19e-4 &  2.33e-4  & 2.33e-4 & 7.60e-6\\ \hline
 94  &      8.8e+4   &  456   &    27.5  &    0.90  &  3.2e-4 & 9.68e-5 & 9.10e-5 & 1.24e-7\\ \hline
 114 &   1.28e+5  &  309   &    42.3  &     1.66 &   1.81e-4&  5.59e-5 & 5.14e-5 & 6.41e-6\\ \hline
 130 &   1.68e+5  &  309  &     69.8 &      2.78&    1.20e-4&   3.74e-5 &  3.40e-5 & 6.98e-7\\ \hline
 145 &   2.08e+5  &  309  &     135.1 &      5.52 & 8.69e-5 & 2.72e-5 & 2.46e-5 &   2.40e-7\\ \hline
 158 &   2.48e+5  &  309  &     296.8  &     7.47 & 6.66e-5 & 2.10e-5 & 1.88e-5 &   8.07e-8\\ \hline
 170 &   2.88e+5  &  309  &     342.7  &     11.18 &  5.32e-5 & 1.68e-5 & 1.50e-5 &   1.80e-8\\ \hline
 182 &   3.28e+5  &  309  &     400.6  &     11.52 & 4.37e-5 & 1.38e-5 & 1.24e-5  &  7.41e-9\\ \hline
 192 &   3.68e+5  &  309  &     417.7  &     13.40 &  3.67e-5 & 1.16e-5  & 1.04e-5  &  4.86e-9\\ \hline
\end{tabular}
\end{center}
\end{table}

As the last experiment we use Hilbert matrix as the
coefficient matrix $A$ in system (\ref{eq:1}), the solution $x$ is exactly
as in example 2.

\begin{figure}[ht]\begin{center}
\includegraphics[height=5cm,width=6cm]{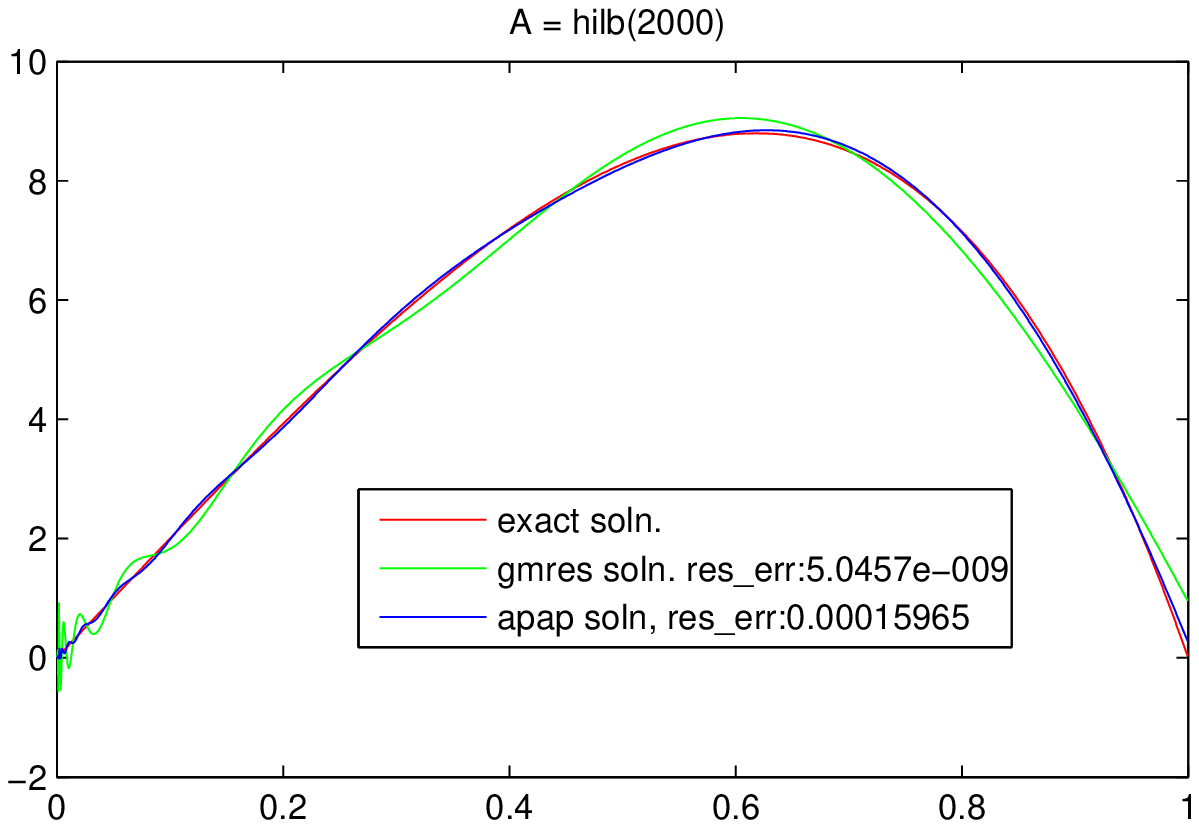}
\includegraphics[height=5cm,width=6cm]{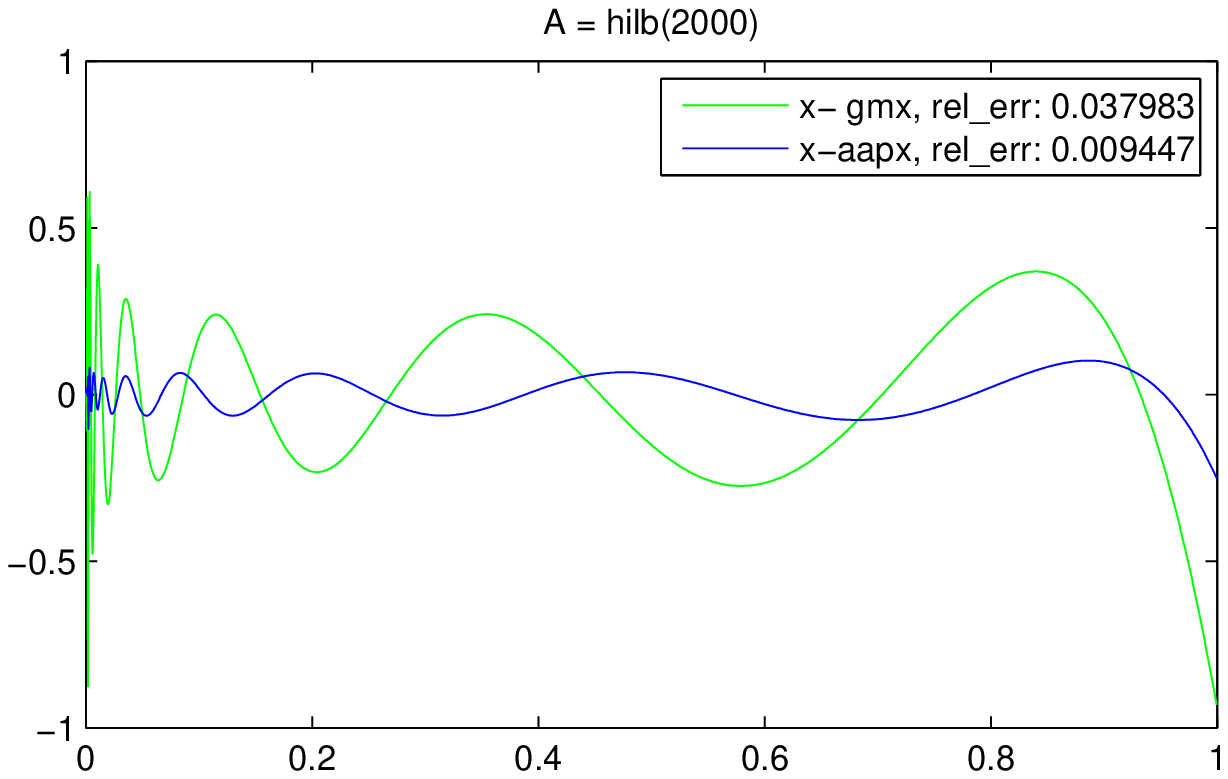}
\caption{Comparison of approx. solns between GMRES and APAP}
\label{fig:hilb2000solns}    
\end{center}\end{figure}

Hilbert matrix is a well-known extremely ill-conditioned matrix and
its condition number grows exponentially. In our experiments  the
direct solver in the MATLAB math package will fail to produce any
significant solution to system (\ref{eq:1}) as $n$ is greater than
$16$. However by using APAP we can solve this system with $n$ up to a
few thousand in this case(see figure \ref{fig:hilb2000solns}).
Again in this case CG can be applied and it takes much less cpu time to reach the same accuracy;

It is interesting to notice here that although the approximate solution given by GMRES yields a much better relative residual, its relative error is a little worse than that of approximate solution given by APAP.

\section{Comments and Summary}

In this paper we discussed a new type of projection methods with the newly introduced AP technique. The major features of these type of iterative methods which make them differ from current existing prevalent Krylov subspace methods includes: (1) the inner products between each approximate solution and the exact solution is
recorded and used for later approximations; (2) they are the first type of non-Krylov subspace methods as far as authors know; (3) the AP techniques actually help to expand the original system into a much larger size of systems (i.e., many more equations can be embedded into the original system with the same solution) and therefore bring much more opportunity for designing accelerative schemes like the one in APAP method.

These type of methods can overcome some shortcomings of current prevailing Krylov subspace methods and exhibit better performance in many of our test problems, especially in case of large sparse linear systems. We have to point out that the construction of some test systems are made so that the exact solutions have dominant components coming from the eigenvalues of the coefficient matrix $A$ with smallest eigenvalues in magnitude, and our test shows that Krylov subspaces usually have a slow convergence speed in these situations, while the APAP method introduced here has a much stable and better performance behavior. APAP can also be used to solve systems with dense coefficient matrices, however to make it applicable, one needs an efficient process to get the projection vectors of $x$ into subspaces formed by row vectors of submatrices of the coefficient matrix, which will be introduced in our later work. When the size of the blocks decreases, or equivalently the number of blocks increases, the convergence speed deteriorate. A remedy is to simply increase the number of AP sweep in each AP process and our numerical experiments show that the time cost is quite reasonable. Currently there is no theoretical results for predicting the iteration numbers needed for any specified tolerance level, since it depends on detailed error analysis of AP process, which seems to be a challenging problem since there does not exist a so-called iteration matrix in the above AP schemes as those in traditional iterative schemes.

\section*{Acknowledgements}
Authors are grateful to the unknown referees for their pertinent suggestions and great help in preparation of this paper.

\bibliographystyle{plain}
\bibliography{peng_refv2}
\end{document}